\documentclass[11pt]{article}
\usepackage{amsthm}
\usepackage{amsfonts}
\usepackage{graphicx}
\newtheorem*{theorem}{Theorem}
\renewcommand{\phi}{\varphi}
\newcommand{\dhat}{\widehat{d}}
\newcommand{\bhat}{\widehat{b}}
\newcommand{\Z}{{\mathbb Z}}
\title{Goldbug Variations}
\author{Michael Kleber \\ The Broad Institute at MIT \\ kleber@broad.mit.edu}
\date{In {\em Mathematical Intelligencer} 27 \#1 (Winter 2005)}
\begin{document}
\maketitle

Jim Propp bugs me sometimes.  I'm usually glad when he does.

Today, Jim's bugs are trained to hop back and forth on the positive
integers: place a bug at $1$, and with each hop, a bug at $i$ moves to
either $i+1$ or $i-2$.  Of course, it might jump off the left edge;
we put two bug-catching cups at $0$ and $-1$.  Once a bug lands in a
cup, we start a new bug at $1$.

\begin{figure}
\includegraphics[width=\textwidth]{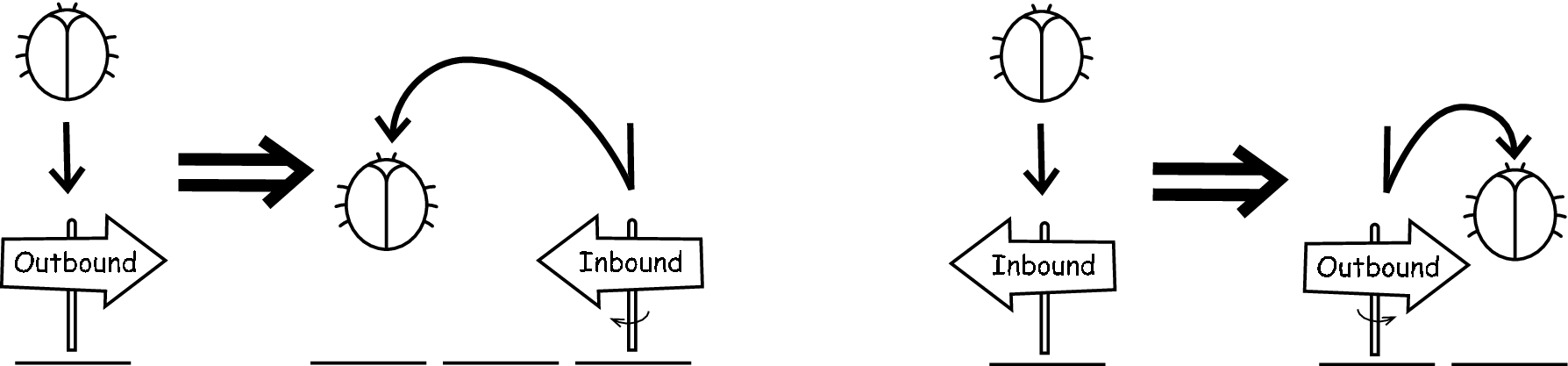}
\caption{The two bug bounces.}
\label{fig:bounces}
\end{figure}

\begin{figure}
\includegraphics[width=\textwidth]{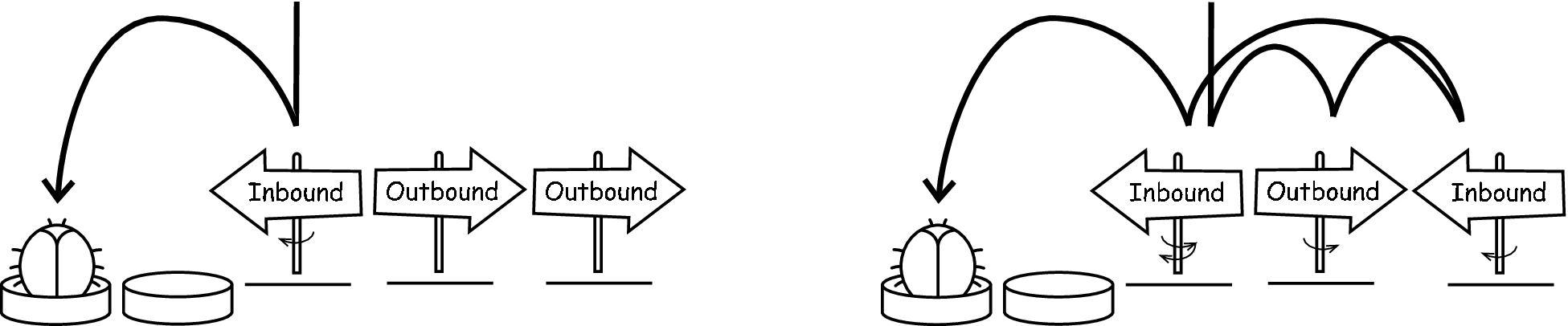}
\caption{The bounces of (a) the first bug, (b) the fourth bug.}
\label{fig:plop}
\end{figure}

What I haven't mentioned is how the bugs decide whether to jump left
or right.  We could declare it a random walk, stepping in either
direction with probability $\!{}^1\!/_2$, but the U.~of Wisconsin professor's 
bugs are more orderly than that.  At each location $i$, there is a signpost 
showing an arrow: it can point either Inbound, towards $i-2$ and the bug 
cups, or Outbound, towards infinity.  The bugs are somewhat contrarian, 
so when a bug lands at $i$, it first {\em flips} the arrow to point the
opposite direction, and then hops that way~(Fig.~\ref{fig:bounces}).  
Add an initial condition that all arrows begin pointing Outbound, and we 
have a deterministic system.  Let bugs hop till they drop~(Fig.~\ref{fig:plop}).

Well, what happens?  Or, for those who would like a more directed
question: First, show that every bug lands in a cup (as opposed to
going off to infinity, or bouncing around in some bounded region
forever).  Second, find what fraction of the bugs end their journeys
in the cup at $-1$, in the many-bugs limit.  Go ahead, I'll wait.  Do the first
ten bugs by hand and look for the beautiful pattern.  You
can even skip ahead to the next section and read about another related
bug, one with far more inscrutable behavior, and come back and read
my solution another day.

I should mention, by way of a delaying tactic, that the analysis of
this bug was first done by Jim Propp and Ander Holroyd.  A previous
and closely related bug of theirs, in which every third visit to $i$
leads to $i-1$, appears in Peter Winkler's new book {\em Mathematical
Puzzles: A Connoisseur's Collection}, yet another delightful
mathematical offering from publishers A~K~Peters \cite{W}.  The book itself
is a gold-mine of the type of puzzles that I expect readers of this
column would enjoy immensely.
Winkler's solutions are insightful, well-written, and often leave the
reader with more to think about than before.  The preceding is an
unpaid endorsement.

Very well, enough filler; here's my answer.  If you solved the problem
without developing something like the theory below, please 
do let me know how.

Before we get to the serious work, let's answer the question, can a
bug hop around in a bounded area forever?  It cannot: let $i$ be
the minimal place visited infinitely often by the bug --- oops, half
the time, that visit is followed by a trip to $i-2$, a contradiction.  
So each bug either lands in a cup or, as far as we know now, runs off 
to infinity.

To motivate what follows, let's have a little inspiration: some
carefully chosen experimental data.  Perhaps you noticed that of the
first five bugs, the cup at $-1$ catches three, while the cup at $0$
catches two.  If that is not sufficiently suggestive, let me mention
that after bug eight the score is $5:3$, while after bug thirteen it
is $8:5$.  On the speculation that this golden ratio trend continues,
let us refer to the bouncing insects as goldbugs.

Now the details.  Suppose that $\phi$ is --- I can only imagine your surprise ---
a real number satisfying $\phi^2=\phi+1$; when I care to be specific
about which root, I will write $\phi_\pm$ for $(1\pm\sqrt5)/2$.

The goldbugs and signs are, in fact, a number written in base $\phi$.
Position $i$ has ``place value'' $\phi^i$, and the digits are
conveniently mnemonic: Outbound is $0$, Inbound is $1$, and the bug
itself is the not entirely standard digit $\phi$, which may appear in
addition to the $0$ or $1$ in the ``$\phi^i$'s place,'' where it contributes
$\phi^{i+1}$ to the total.  Of course, numbers do not have
a unique representation with this setup, but that's quite deliberate:
the total value is an invariant, unchanged by bug bounces.
\begin{itemize}
\item{Bounce left:} Suppose the bug arrives in position $i$ and there
  is an Outbound arrow.  The value of this part of the configuration
  is $(\phi^i\times\phi)+(\phi^i\times0)$.  After the bounce, position $i$
  holds an Inbound arrow and the bug has bounced to position $i-2$,
  for a total value of $(\phi^i\times1)+(\phi^{i-2}\times\phi)$.  And these are the
  same, since $\phi^{i-1}+\phi^i=\phi^{i+1}$.
\item{Bounce right:} Suppose the bug arrives in position $i$ and there
  is an Inbound arrow.  The value before the bounce is
  $(\phi^i\times\phi)+(\phi^i\times1)$.  After the bounce, the arrow points Outbound,
  value zero, and the bug is at $i+1$, for a value of $\phi^{i+2}$,
  again unchanged.
\end{itemize}

Now let's see what happens when we add a bug at $1$, let it run its
course, and then remove it after it lands in a bug cup.
Placing a bug at $1$ increases the system's value by $\phi^2$.  
If it eventually lands in the cup at $0$ and is then removed, the value 
of the system drops by $\phi$, while if it lands in the cup at $-1$ and is
removed, the value of the system drops by $\phi^0=1$.  So the {\em net} 
effect of adding a bug at 1, running the system, and then removing the bug 
is that the value of the system increases by $\phi^2-\phi=1$ for cup $0$, 
and by $\phi^2-1=\phi$ for cup $-1$.

Now we are well-equipped to answer the original questions, and then some.

\begin{itemize}
\item
Can a bug run off to infinity?  It cannot: if we take $\phi$ to mean
$\phi_+ \approx 1.61803$, then each bug can increase the net value of the
system by at most $\phi_+$, and the positions far to the right are
inaccessible to the goldbugs, because they would make the value of the 
system too high.  So every bug lands in a cup.

\item
What is the ratio of bugs landing in the two cups?  This time for $\phi$, think
about $\phi_- \approx -.61803$.  Between bugs, when the system 
consists of just Inbound and Outbound flags, its minimum possible value 
is $\phi_-^1+\phi_-^3+\phi_-^5+\cdots = -1$, while its maximum possible
value is $\phi_-^2+\phi_-^4+\phi_-^6+\cdots=-\phi_-\approx .61803$,
which I will write instead as $1/\phi_+$ to help us remember its sign.

So the value of the system is trapped in the interval $[-1,1/\phi_+]$,
and as each successive bug passes through the system, the value
either increases by $1$, or else increases 
by $\phi_-$, i.e.\ decreases by $1/\phi_+$.  
For the value to stay in any bounded region, it must, in
the limit, step down $\phi_+$ times as often as it steps up, and so the
bugs land in the $-1$ cup $1/\phi_+$ of the time.  Moreover, this
approximation is very good: if $a$ bugs have landed at $-1$ and $b$
have landed at $0$, the system's value $b-a/\phi_+$ must lie in our
interval, so $|b/a - 1/\phi_+| < 1/a$.  Every single approximation
$b/a$ is one of the two best possible given the denominator, and that
denominator grows as $n/\phi_+$.

In fact, notice that the length of the interval $[-1,1/\phi_+]$ is the
sum of the two jump sizes, the smallest we could possibly hope for.
If the value lies in the bottom $1/\phi_+$ of the interval, it must
increase by $1$, while if it lies in the top $1$ of the interval, it
must decrease by $1/\phi_+$.  This leaves a single point, $1/\phi_+-1
\approx -.38297$, where the bug's destination cup isn't determined.
But that value can be attained only with infinitely many Inbound
signs: if we run a single bug through a system with that starting
value, its ending value would be either $1/\phi_+$ or $-1$ --- and to
attain those two bounds, we found above, you must sum the infinite
series of all positive or negative place values.  So for any initial
configuration with only finitely many signs pointing Inbound, the
configuration's numerical value alone determines the destination cup of 
every single bug; you don't need to keep track of the arrows at all.

\item
If you prefer integers to irrationals, consider instead the following
invariant.  Label the (bug cups and the) sites with the Fibonacci 
numbers (0,1,)1,2,3,5,8,$\ldots$, as in Figure~\ref{fig_fibonacci}.
Given an Inbound arrow the value $F_i$ of its site, and give the
bug there the value $F_{i+1}$ of the site to its right.  This invariant, 
of course, is an appropriate linear combination of the $\phi_+$ and 
$\phi_-$ ones above.

Adding a bug to the system now increases the value by 2, but what's
special is that removing the bug at either $0$ {\em or} $-1$ subtracts 1.  
So the bugs implement an accumulator: after the $n$th bug lands in 
its cup, we can look at the signs it has left behind, and read off the 
number $n$, written in base Fibonacci!

\begin{figure}
$$
\newcommand{\I}{\stackrel{\rm I}\longleftarrow}
\renewcommand{\O}{\stackrel{\rm O}\longrightarrow}
\begin{array}{ccccccccccc}
\cup & \cup & \O & \I & \O & \I & \I & \I & \O & \I & \I \\
{\scriptscriptstyle (0\,} & {\scriptscriptstyle \,1)} & 1 & 2 & 3 & 5 & 8 & 13 & 21 & 34 & 55
\end{array}
$$
\caption{The signs after bug 117 passes through the system;
117=2+5+8+13+34+55.  Representations in base Fibonacci
are not unique; ours is characterized by a lack of two consecutive
zeros.}
\label{fig_fibonacci}
\end{figure}

\ldots Hold the presses!  Matthew Cook has pointed out to me that this
point of view can be taken further.  We still get an invariant if we shift
all those Fibonacci labels one site to the right.  Now adding a bug
increases the system's total value by 1, and removing it from the right 
bug cup decreases the value by 1 --- but removing it from the left 
bug cup decreases the score by 0.  So after $n$ bugs pass through,
the value of the Inbound arrows they leave behind counts the number
of bugs that ended their journeys in the left bug cup.

Furthermore we can shift the labels right a second time.  Now when
the bug lands in the left cup, its value must be the $-1$st Fibonacci
number, before $0$ --- which is 1 again, if the Fibonacci recurrence
relation is still to hold.  So with these labels, the Inbound arrows will
count the number of bugs which landed in the right cup.  As an
exercise, decode the arrows to learn that the 117 bugs leading to 
Figure~\ref{fig_fibonacci} split 72:45.

Once we know that the same set of arrows simultaneously count the 
total, left-cup, and right-cup bugs, it's straightforward to see that
the ratios of these three quantities are the same as the ratios of three
consecutive Fibonacci numbers, in the limit.

\end{itemize}

These arrow-directed goldbugs are doing a great job of what what Jim Propp 
calls ``derandomization.''  It's straightforward to analyze the corresponding 
random walk, in which bugs hop left or right each with probability $\!{}^1\!/_2$:
Let $p_i$ be the probability that a bug at place $i$ eventually ends up in
the cup at $-1$, and solve the recurrence $p_i=(p_{i-1}+p_{i+1})/2$ with
$p_{-1}=1$, $p_0=0$ --- oh, and make up for one too few initial conditions
by remembering that all the $p_i$ are probabilities, so no larger than 1.
Here too we get $p_1=1/\phi_+$.  

So the deterministic goldbugs have the same limiting behavior as their 
random-walking cousins.  But if we run the experiment with $n$ random
walkers and count the number of bugs in a cup, we'll generally see 
variation on the order of $\sqrt{n}$ around the expected number.
Remember the sharp results of the $\phi_-$ invariant: $n$ goldbugs, 
by contrast, simulate the expected behavior with only constant error!

There are more results on these and related one-dimensional 
not-very-random walks.  But let's move on --- these bugs long for
some higher-dimensional elbow room.

\section*{The Rotor-Router}

This time, we will set up a system with bugs moving around the integer
points in the plane.  It will differ from the above in that this will
be an {\em aggregation model}.  Bugs still get added repeatedly at one
source, but instead of falling into sinks (bug cups), the bugs will
walk around until they find an empty lattice point, then settle down
and live there forever.

\begin{figure}
\includegraphics[width=\textwidth]{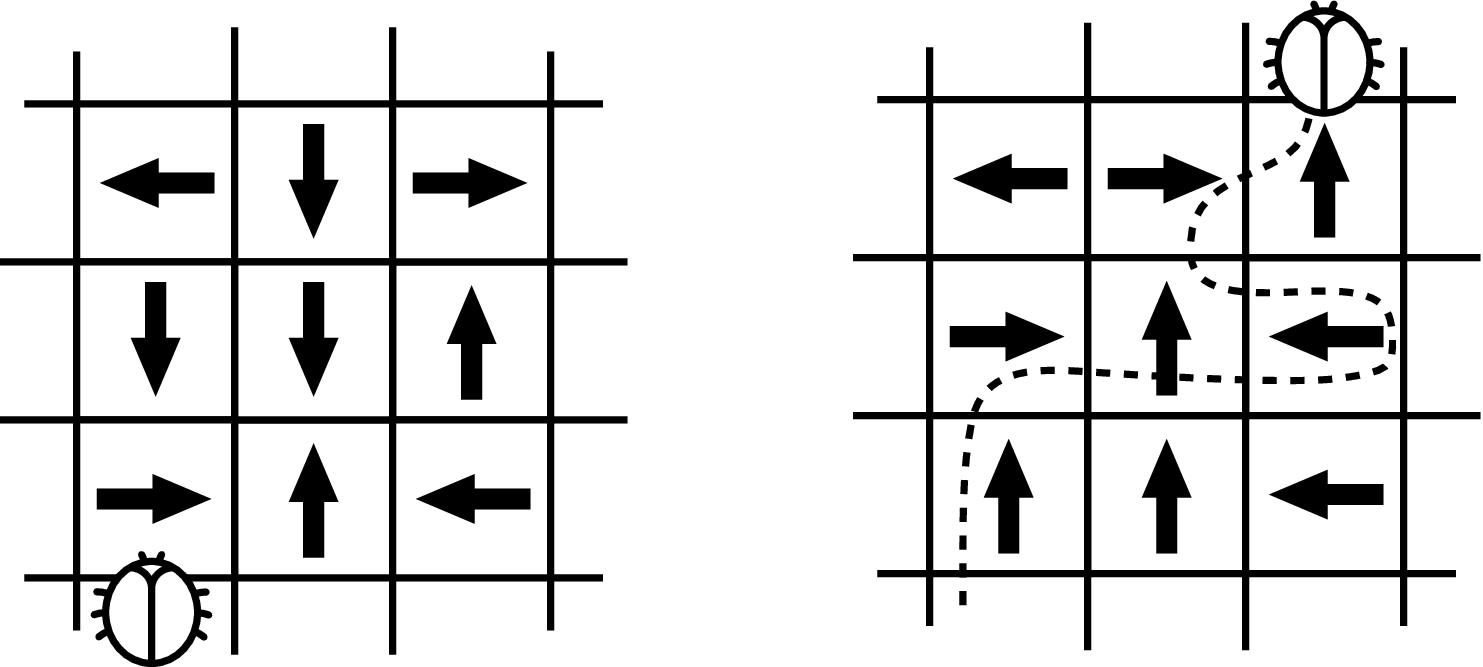}
\caption{A bug's ramble through some rotors: before and after.}
\label{fig:ramble}
\end{figure}

Generalizing the Inbound and Outbound arrows of the goldbugs, we
decree that each lattice site is equipped with an arrow, or rotor,
which can be rotated so that it points at any one of the four 
neighbors.  (Propp uses the word rotor to refer to the
two-state arrows in dimension 1 as well.)
The first bug to arrive at a particular site occupies it forever, and
we decree that it sets the arrow there pointing to the East.  Any bug
arriving at an occupied site rotates the arrow one quarter turn
counterclockwise, and then moves to the neighbor at which the rotor
now  points --- where it may find an empty site to inhabit, or may find a new 
arrow directing its next step.  In short, the first bug to reach $(i,j)$ lives there,
and bugs which arrive thereafter are routed by the rotor: first to the
North, then West, South, and East, in that order, and then the cycle
begins with North again~(Fig~\ref{fig:ramble}).

Once a bug finds an empty site to inhabit, we drop a new bug at the
origin, and this one too meanders through the field of rotors, both
directed by the arrows and changing them as a result of its visits.
Every bug does indeed find a home eventually, and the proof is the
same as for the goldbugs: the set of all sites which a particular bug
visits infinitely often cannot have any boundary.

And so we ask, as we add more and more bugs, what does the set of
occupied sites look like?

\begin{figure}
\vspace{-1.5in}
\begin{center}
\makebox[0pt]{\includegraphics[width=8in]{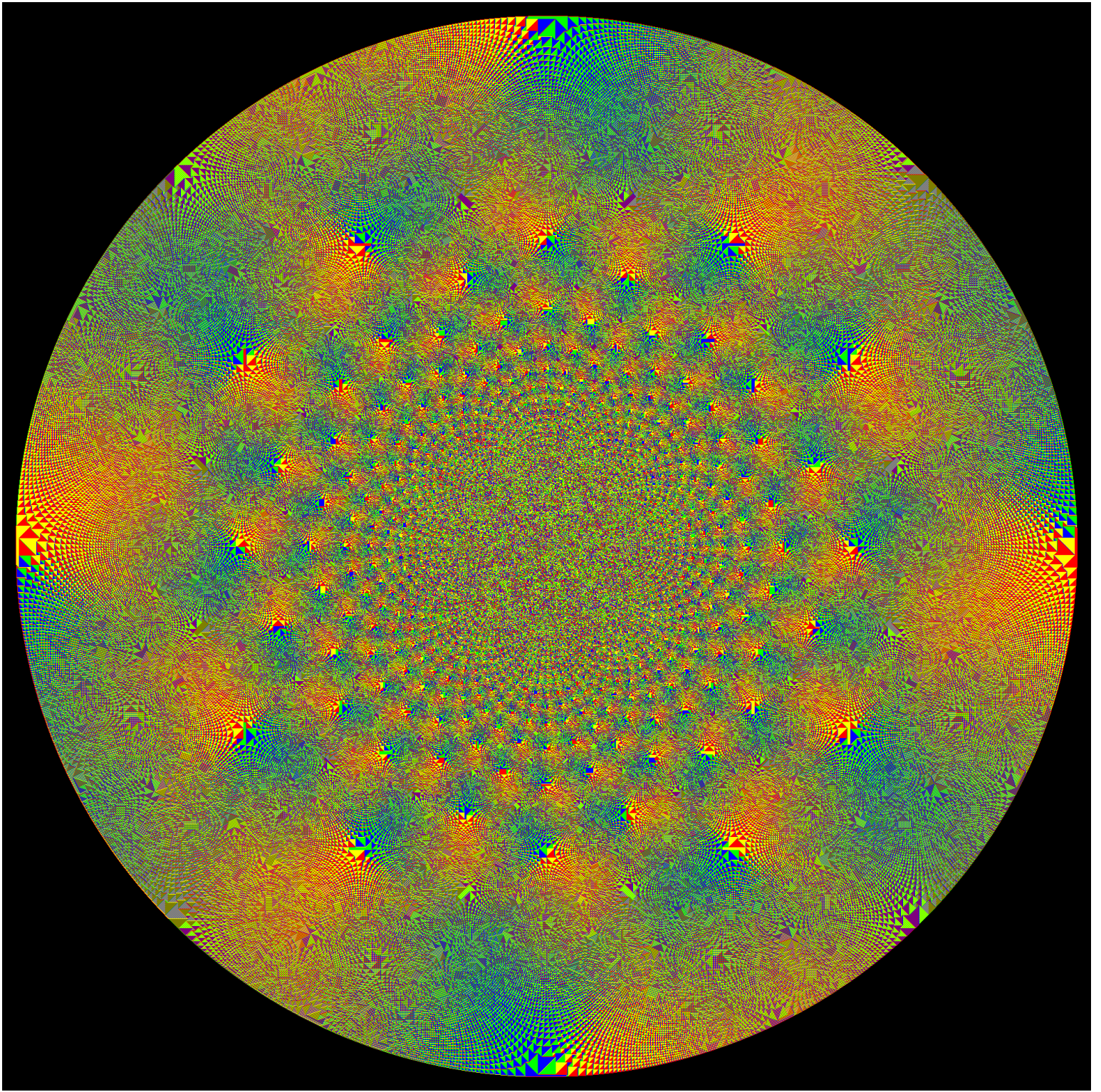}}
\end{center}
\caption{The rotor-router blob after three million bugs.  If you are viewing
this image on-line, try magnifying to the point where you can see individual
pixels as small boxes.  Also, for best appearance of the entire image, use
a viewer that will smooth/dither/antialias graphics.}
\label{fig:3million}
\end{figure}

Let's take a look at the experimental answer.  The beautiful image you
see in Figure~\ref{fig:3million} is a picture of the set of occupied sites after three million bugs have 
found their permanent homes.  The sites in black are vacant, still awaiting 
their first visitor.  Other sites are colored according to the direction of their
rotor, red/yellow for East/West and green/blue for North/South.
On the cover \marginpar{\tiny\raggedright\sloppypar\em Close-up not included 
in the web version, but you should be able to zoom in on Fig.~\ref{fig:3million} for
the same effect.} is a close-up of a part of the boundary of the occupied region.
At \verb"http://www.math.wisc.edu/~propp/rotor-router-1.0/" you can find
a Java applet by Hal Canary and Francis Wong, if you'd like to experiment
yourself.

As you can see, the edge of the occupied region is extraordinarily 
round: with three million bugs, the occupied site furthest from the origin is 
at distance $\sqrt{956609}$, while the unoccupied site nearest the origin is at
distance $\sqrt{953461}$, a difference of about 1.6106.

And the internal coloration puts on a spectacular display of both large-scale 
structure and intricate local patterns.  When Ed Pegg featured the rotor-router on 
his Mathpuzzle web site, he dubbed the picture a Propp Circle, and to this day I
am jealous that I didn't think of the name first: it connotes precisely the right
mixture of aesthetic appreciation and conviction that there must be something
deep and not fully understood at work.

Recall, for emphasis, that this was formed by bugs walking deterministically on 
a square lattice, not a medium known for growing perfect circles.  Moreover, the 
rule that governed the bugs' movement is inherently asymmetric: every site's
rotor begins pointing East, so there was no guarantee that the set of
occupied sites would even appear symmetric under rotation by $90^\circ$,
much less by unfriendly angles.  On the other hand, while the overall shape 
seems to have essentially forgotten the underlying lattice, the internal structure 
revealed by the color-coded rotors clearly remembers it.

Lionel Levine, now a graduate student at U.C.~Berkeley, wrote an undergraduate
thesis with Propp on the rotor-router and related models~\cite{levine}.  It contains the best 
result so far on the roundness of the rotor-router blob: after $n$ bugs, it
contains a disk whose radius grows as $n^{1/4}$.
Below I report on two remarkable theorems which do not quite prove that the 
rotor-router blob is round, but which at least make me feel that it ought to be.
I have less to offer on the intricate internal structure, but there is a connection
to something a bit better understood.  Let's get to work.

\subsection*{IDLA}

Internal Diffusion Limited Aggregation is the random walk--based
model which Propp de-randomized to get the rotor-router.  Most 
everything is as above --- the plane starts empty, add bugs to the 
origin one at a time, each bug occupies the first empty site it reaches.  
But in IDLA, a bug at an occupied site walks to a random neighbor, each 
with probability $\!{}^1\!/_4$.

The idea underlying IDLA comes from a paper by Diaconis and Fulton~\cite{DF}.
They define a ring structure on the vector space whose basis is labeled by
the finite subsets of a set $X$ equipped with a
random walk.  To calculate the product of subsets $A$ and $B$, begin with the set
$A\cup B$, place bugs at each point of $A\cap B$, and allow each bug to execute a
random walk until it reaches an outside point, which is then added to the set.
The product of $A$ and $B$ records the probability distribution of possible
outcomes.  This appears to depend on the order in which the bugs do their 
random walks, but in fact it does not --- we'll explore this theme soon.

Consider the special case where $X$ is the $d$-dimensional integer
lattice with the random walk choosing uniformly from among the nearest
neighbors.  Then repeatedly multiplying the singleton $\{{\mathbf
0}\}$ by itself is precisely the random-walk version of the
rotor-router.  A paper of Lawler, Bramson and Griffeath~\cite{LBG}
dubbed this Internal Diffusion Limited Aggregation, to emphasize
similarity with the widely-studied Diffusion Limited Aggregation model
of Witten and Sandler~\cite{DLA}.  DLA simulates, for example, the
growth of dust: successive particles wander in ``from infinity'' and
stick when they reach a central growing blob.  The resulting growths
appear dendritic and fractal-like, but rigorous results are hard to
come by.

In contrast, the growth behavior of IDLA has been rigorously
established.  It is intuitive that the growing blob should be
generally disk-shaped, since the next particle is more likely to fill
in an unoccupied site close to the origin than one further away.  But
the precise statement in~\cite{LBG} is still striking: the random walk
manages to forget the anisotropy of the underlying lattice entirely!
\begin{theorem}[Lawler--Bramson--Griffeath]
Let $\omega_d$ be the volume of the $d$-dimensional unit sphere.
Given any $\epsilon>0$, it is true with probability 1 that for all
sufficiently large $n$, the $d$-dimensional IDLA blob of $\omega_d
n^d$ particles will contain every point in a ball of radius
$(1-\epsilon)n$, and no point outside of a ball of radius
$(1+\epsilon)n$.
\end{theorem}

To be more specific, we could hope to define inner and outer error
terms such that, with probability 1, the blob lies between the balls
of radius $n-\delta_I(n)$ and $n+\delta_O(n)$.  In a subsequent
paper~\cite{L}, Lawler proved that these could be taken on the order
of $n^{1/3}$.  Most recently, Blach\`ere~\cite{B} used an induction argument
based on Lawler's proof to show that these error terms were even
smaller, of logarithmic size.  The precise form of the bound changes
with dimension; when $d=2$ he shows that $\delta_I(n) = O((\ln n \ln
(\ln n))^{1/2})$ and $\delta_O(n) = O((\ln n)^2)$.  Errors on that order 
were observed experimentally by Moore and Machta~\cite{MM}.

So how does the random walk--based IDLA relate to the deterministic
rotor-router?  We start drawing the connection with one key fact.

\subsection*{It's abelian!}

Here's a possibly unexpected property of the rotor-router model: it's
abelian.  There are several senses in which this is true.

Most simply, take a state of the rotor-router system --- a set of occupied
sites and the directions all the rotors point --- and add one bug at a
point $P_0$ (not necessarily the origin now) and let it run around and 
find its home $P_1$.  Then add another at $Q_0$ and let it run until it 
stops at $Q_1$.  The end state is the same as the result of adding the two 
bugs in the opposite order.  

This relies on the fact that the bugs are indistinguishable.  Consider the 
(next-to-)simplest case, in which the paths of the $P$ and $Q$ bugs cross at
exactly one point, $R$.  If bug $Q$ goes first instead, it travels from $Q_0$ to 
$R$, and then follows the path the $P$ bug would have, from $R$ to $P_1$. 
The $P$ bug then goes from $P_0$ to $R$ to $Q_1$.  At the place where their
paths would first cross, the bugs effectively switch identities.  For paths whose
intersections are more complicated, we need to do a bit more work, but the 
basic idea carries us through.

Taking this to an extreme, consider the ``rotor-router swarm''
variant, where traffic is still directed by rotors at each lattice
site, but any number of bugs can pass through a site simultaneously.
The system evolves by choosing any one bug at random and moving it one
step, following the usual rotor rule.  Here too the final state is
independent of the order in which bugs move; read on for a proof.  
To create our rotor-router picture, we can place three million bugs at the origin 
simultaneously, and let them move one step at a time, following the rotors, 
in whatever order they like.

In fact, even strictly following the rotors is unnecessary.  The
rotors control the order in which the bugs depart for the various
neighbors, but in the end, we only care about how many bugs head in
each direction.

Imagine the following setup: we run the original rotor-router
with three million bugs as first described, but each time a bug
leaves a site, it drops a card there which reads ``I went North'' or
whichever direction.  Now forget about the bugs, and look only at the
collection of cards that were left behind at each site.  This certainly
determines the final state of the system: a site ends up occupied if
and only if one of its neighbors has a card pointing towards it.

Now we could re-run the system with no rotors at all.  When a bug
needs to move on, it may pick up {\em any} card from the site it's on
and move in the indicated direction, eating the card in the process.
No bug can ever ``get stuck'' by arriving at an occupied site with no 
card to tell it a way to leave: the stack of cards at a given site is just
the right size to take care of all the bugs that can possibly arrive there
coming from all of the neighbors.  (There is, however, no guarantee
that all the cards will get used; left-overs must form loops.)  
A version of this ``stacks of cards'' idea appeared in Diaconis and 
Fulton's original paper, in the proof that the random walk version is 
likewise abelian --- i.e.\ that their product operation is well-defined.

If the bugs are so polite as to take the cards in the cyclic N-W-S-E order 
in which they were dropped, then we simulate the rotor-router exactly.
If we start all the bugs at the origin at once and let them move in whatever
order they want --- but insist that they always use the top card from the
site's stack --- we get the rotor-router swarm variant above; QED.

%If instead the bugs choose randomly from among the remaining cards, we get
%something like IDLA, but with the random walk subject to the Gambler's 
%Fallacy.  In a casino, I wouldn't recommend the logic ``This roulette wheel just 
%came up Black ten  times in a row; time to bet big on Red!''  But as cards get 
%used up, a bug choosing randomly from what remains is more likely to leave 
%a site headed in a direction underrepresented so far.

\subsection*{Rotor-roundness}

Now let me outline a heuristic argument that the rotor-router blob ought
to be round, letting the Lawler--Bramson--Grifffeath paper do all the heavy
lifting.

I'd like to say that, for any $c<1$, the $n$-bug rotor-router blob contains 
every lattice site in the disk of area $cn$ --- as long as $n$ is sufficiently 
large.  My strategy is easy to describe.  Just as we did four paragraphs ago,
think of each lattice site as holding a giant stack of cards: one card
for each time a bug departed that site while the $n$-bug rotor-router blob grew.
Now we start running IDLA: we add bugs at the origin, one at a time,
and let them execute their random walk.  But each time a bug randomly
decides to step in a given direction, it must first look through the
stack of cards at its site, find a card with that direction written on
it, and destroy it.

As long as the randomly-walking bugs always find the cards they look
for, the IDLA blob that they generate must be a subset of the
rotor-router blob whose growth is recorded in the stacks of cards.
This key fact follows directly from the abelian nature of the models.

So the central question is, how long will this IDLA get to run before
a bug wants to step in a particular direction and finds that there is
no corresponding card available?  Philosophically, we expect the
IDLA to run through ``almost all'' the bugs without hitting such a
snag: for any $c<1$, we expect $cn$ bugs to aggregate, as long as $n$
is sufficiently large.  If we can show this, we are certainly done:
the rotor-router blob contains an IDLA blob of nearly the same area,
which in turn contains a disk of nearly the same area, with
probability one.

To justify this intuition, we clearly need to examine the function $d(i,j)$ which
counts the number of departures from each site.  This is a nonnegative
integer--valued function on the lattice which is almost harmonic, away
from the origin: the number of departures from a given site is about
one quarter of the total number of visits to its four neighbors.
$$
d(i,j) \approx \frac{d(i+1,j)+d(i-1,j)+d(i,j+1)+d(i,j-1)}{4} - b_{(i,j)}
$$
Here $b_{(i,j)}=1$ if $(i,j)$ is occupied and $0$ otherwise, to account
for the site's first bug, which arrives but never departs.
When $(i,j)$ is the origin, of course, the right-hand side should be
increased by the number of bugs dropped into the system.  
Matthew Cook calls this the ``tent equation'': each site is forced to be a 
little lower than the average of its neighbors, like the heavy fabric of a 
circus tent; it's all held up by the circus pole at the origin --- or perhaps by
a bundle of helium balloons which can each lift one unit of tent
fabric, since we do not get to specify the height of the origin, but rather
how much higher it is than its surroundings.
 
For the rotor-router, the approximation sign above hides some rounding 
error, the precise details of which encapsulate the rotor-router rule.  For
IDLA, this is exact if we replace $d$ by $\dhat$, the {\em expected}
number of departures, and replace $b$ by $\bhat$, the probability
that a given site ends up occupied.  (The results of~\cite{LBG} even give
an approximation of $\dhat$.)

Now, I'd like to say that at any particular site, the mean number of
departures for an IDLA of $cn$ bugs (for any $c<1$ and large $n$) should 
be less than the actual number of departures for a rotor-router of $n$ bugs.  
If so, we'd be nearly done, with a just a bit of easy calculation to show that the 
$\sqrt{d}$-sized error terms at each site in the random walk is thoroughly 
swamped by the $(1-c)n$ extra bugs in the rotor-router.  

But this begs the question of showing that the rotor-router's function $d$ and 
the IDLA's function $\dhat$ are really the same general shape.  Their difference 
is an everywhere {\em almost}-harmonic function with zero at the boundary --- but 
to paraphrase Mark Twain, the difference between a harmonic function and an
almost-harmonic function is the difference between lightning and a lightning bug.

\subsection*{Simulation with constant error}

After I wrote the preceding section, I learned of a brand-new result of 
Joshua Cooper and Joel Spencer.  It doesn't turn my hand-waving into a
genuine proof, but it gives me hope that doing so is within reach.
Their paper \cite{CS} contains an amazing result on the relationship between
a random walk and a rotor-router walk in the $d$-dimensional integer lattice $\Z^d$.

Generalizing the rotor-router bugs above, consider a lattice $\Z^d$ in which each 
point is equipped with a rotor --- that is to say, an arrow which points towards one of 
the $2d$ neighboring points, and which can be incremented repeatedly, causing it to 
point to all $2d$ neighbors in some fixed cyclic order.  The initial states of the rotors 
can be set arbitrarily.

Now distribute some finite number of bugs arbitrarily on the points.  We can let this
distribution evolve with the bugs following the rotors: one step of evolution consists
of every bug incrementing and then following the rotor at the point it is on.  (Our 
previous bugs were content to stay put if they were at an uninhabited site, but in this 
version, every bug moves on.)  Given any initial distribution of bugs and any initial 
configuration of the rotors, we can now talk about the result of $n$ steps of 
rotor-based evolution.

On the other hand, given the same initial distribution of bugs, we could just as well
allow each bug to take an $n$-step random walk, with no rotors to influence its
movement.  If you believe my heuristic babbling above, then it is reasonable to
hope that $n$ steps of rotor evolution and $n$ steps of random walk would lead
to similar ending distributions.  

With one further assumption, this turns out to be true in the strongest of senses.
Call a distribution of bugs ``checkered'' if all bugs are on vertices of the same
parity --- that is, the bugs would all be on matching squares if $\Z^d$
were colored like a checkerboard.

\begin{theorem}[Cooper--Spencer]
There is an absolute constant bounding the divergence between the rotor
and random walk evolution of checkered distributions in $\Z^d$, depending 
only on the dimension $d$.
That is, given any checkered initial distribution of a finite number of bugs in $\Z^d$, 
the difference between the actual number of bugs at a point $p$ after $n$ 
steps of rotor-based evolution, and the expected number of bugs at $p$ after
an $n$-step random walk, is bounded by a constant.  This constant is
independent of the number of steps $n$, the initial states of the rotors, and
the initial distribution of bugs!
\end{theorem}

I am enchanted by the reach of this result, and at the same time intrigued
by the subtle ``checkered'' hypothesis on distributions.  (Not only initial
distributions: since each bug changes parity at each time step, a 
configuration can never escape its checkered past.)  The authors
tell me that without this assumption, one can cleverly arrange squadrons of
off-parity bugs to reorient the rotors and steer things away from random
walk simulation.

Thus the rotor-router deterministically simulates a random walk process with
constant error --- better than a single instance of the random process usually
does in simulating the average behavior.  Recall that we saw a similar outcome
in one dimension, with the goldbugs.

There are other results which likewise demonstrate that derandomizing systems
can reduce the error.  Lionel Levine's thesis~\cite{levine} analyzed a type of
one-dimensional derandomized aggregation model, and showed that it can
compute quadratic irrationals with constant error, again improving on the 
$\sqrt{n}$-sized error of random trials.  Joel Spencer tells me that he can use
another sort of derandomized one-dimensional system to generate binomial 
distributions with errors of size $\ln n$ instead of $\sqrt{n}$.  Surely the 
rotor-router should be able to cut IDLA's already logarithmic-sized variations
down to constant ones.  Right?

\section*{Coda: Sandpiles}

All of the preceding discussion addresses the overall shape of the rotor-router
blob, but says nothing at all about the compelling internal structure that's visible
when we four-color the points according to the directions of the rotors.  When we
introduced the function $d(i,j)$, counting the number of departures from the $(i,j)$
lattice site, we were concerned with its approximate large-scale shape, which
exhibits some sort of radial symmetry.  The direction of the rotor tells you the value
of $d(i,j) \bmod 4$, and the symmetry of these least significant bits of $d$ is an
entirely new surprise.

I can't even begin to explain the fine structure --- if you can, please let me know!
But I can point out a surprising connection to another discrete dynamical system,
also with pretty pictures.

Consider once again the integer points in the plane.  Each point now holds
a pile of sand.  There's not much room, so if any pile has five or more grains
of sand, it collapses, with four grains sliding off of it and getting dumped on
the point's four neighbors.  This may, in turn, make some neighboring piles
unstable and cause further topplings, and so on, until each pile has size at
most four.

\begin{figure}
\vspace{-1.5in}
\begin{center}
\makebox[0pt]{\includegraphics[width=8in]{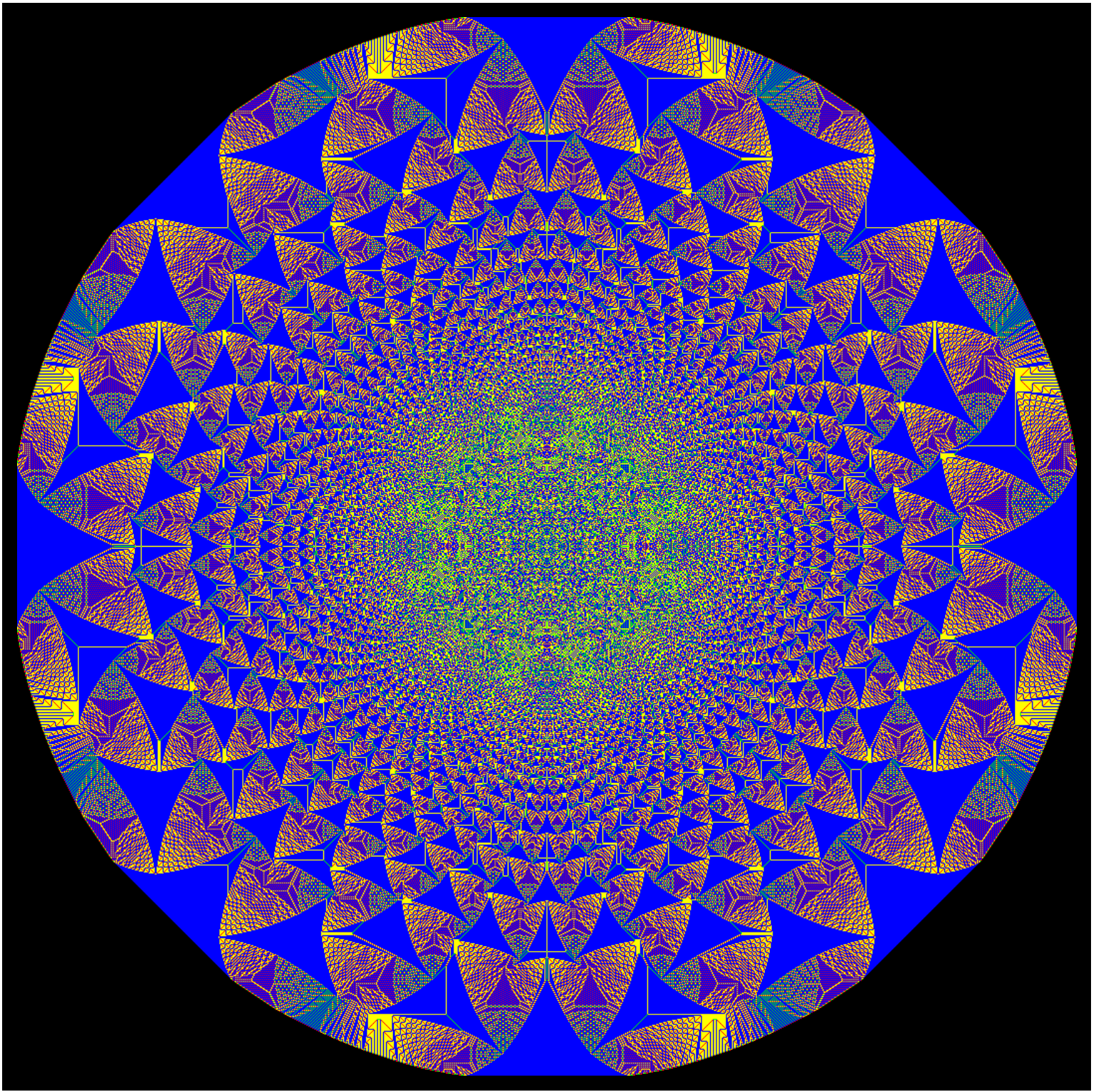}}
\end{center}
\caption{The greedy sand-pile with three million grains.}
\label{fig_gsp}
\vspace{2in}
\end{figure}

Our question: what happens if you put, say, a million grains of sand at the 
origin, and wait for the resulting avalanche to stop?  I won't keep you
hanging; a picture of the resulting rubble appears as Figure~\ref{fig_gsp}.  Pixels
are colored according to the number of grains of sand there in the final
configuration.  The dominant blue color corresponding to the largest
stable pile, four grains.  (This makes some sense, as the interior of such
a region is stable, with each site both gaining and losing four grains,
while evolution happens around the edges.)

This type of evolving system now goes under the names ``chip-firing model''
and ``abelian sandpile model;'' the adjective abelian is earned because the 
operations of collapsing the piles at two different sites commute.
In full generality, this can take place on an arbitrary graph, with an excessively 
large sandpile giving any number of grains of sand to each of its neighbors, 
and some grains possibly disappearing permanently from the system.  Variations
have been investigated by combinatorists since about 1991~\cite{BLS};
they adopted it from the mathematical physics community, which had been
developing versions since around 1987~\cite{BTW,Dhar}.  This too was a
rediscovery, as it seems that the mechanism was first described, under
the name ``the probabilistic abacus,'' by Arthur Engel in 1975 in a math
education journal~\cite{E1,E2}.

I couldn't hope to survey the current state of this field here, or even give proper
references.
The bulk of the work appears to be on what I think of as steady-state questions, 
far from the effects of initial conditions:
point-to-point correlation functions, the distribution of sizes of avalanches, 
or a marvelous abelian group structure on a certain set of recurrent 
configurations.

\begin{figure}
\vspace{-1.5in}
\begin{center}
\makebox[0pt]{\includegraphics[width=8in]{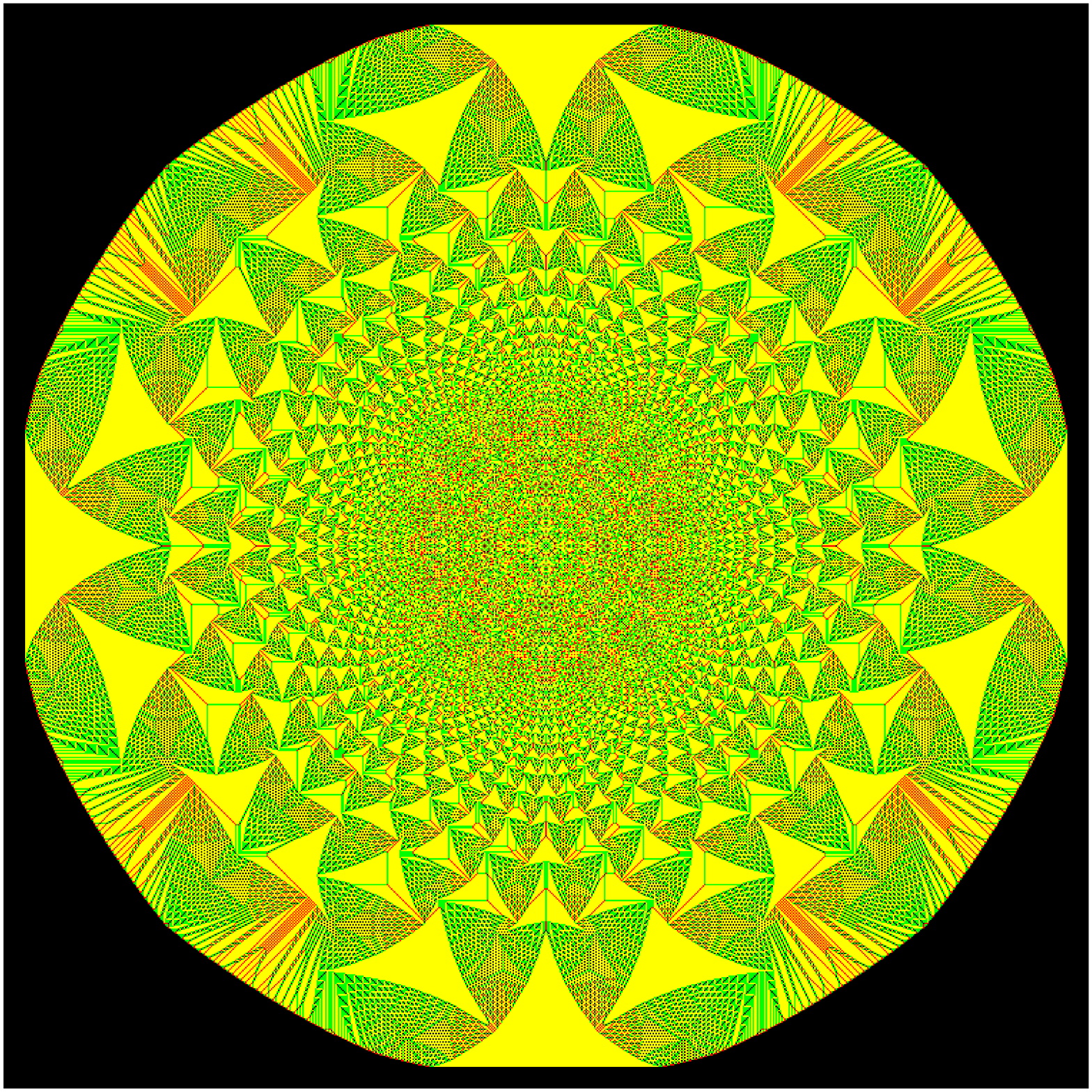}}
\end{center}
\caption{A non-greedy sand-pile.  Here the dominant color is yellow, which
again indicates the maximal stable site, now with three grains.  It is hard to see
the interior pixels colored black, indicating sites which were once filled but are
now empty, impossible in the greedy version.}
\label{fig_sp}
\end{figure}

Our question seems to have a different flavor.  For example, in most
sandpile work, one can assume without loss of generality that a pile collapses
as soon as it has enough grains of sand to give its neighbors what they are owed,
leaving itself vacant.  The version I described above is what I'll call a ``greedy 
sandpile,'' in which each site hoards its first grain of sand, never letting it go.  
The shape of the rubble in Figure~\ref{fig_gsp} does depend on this detail; Figure~\ref{fig_sp}
is  the analogue where a pile collapses as soon as it has four grains, leaving itself 
empty.

Most compelling to me is the fine structure of the sandpile picture.  I'm amazed
by the appearance of fractalish self-similarity at different scales despite the
single-scale evolution rule; I think this is related to what the mathematical 
physics people call ``self-organizing criticality,'' about which I know nothing at all.
But personally, in both pictures I am drawn to the eight-petalled
central rosette, the boundary of some sort of phase change in their internal 
structures.

\subsection*{Bugs in the sand}

So what is the connection between the greedy sandpile and the rotor-router?
Recall the swarm variant of rotor-router evolution: we can place all the bugs at
the origin simultaneously, and let them take steps following the rotor rule in any 
order, and still get the same final state.  

Since we get to choose the order, what if we repeatedly pick a site with at least
four bugs waiting to move on, and tell four of them to take one step each?  
Regardless of its state, the rotor directs one to each neighbor, and we mimic 
the evolution rule of the greedy sandpile perfectly.  If we keep doing this until
no such sites remain, we realize the sandpile final state as one step along one
path to the rotor-router blob.

Note, in particular, that the $n$-bug rotor-router blob must contain all sites in the
$n$-grain greedy sandpile.  Surely it should therefore be possible to show that both
contain a disk whose radius grows as $\sqrt{n}$.

More emphatically, the sandpile performs precisely that part of the evolution of the 
rotor-router that can take place without asking the rotors to break symmetry.  If we
define an energy function which is large when multiple bugs share a site, then the
sandpile is the lowest-energy state which the rotor-router can get to in a completely
symmetric way.  

When we invoke the rotors, we can get to a state with minimal energy but without 
the a priori symmetry that the sandpile evolution rule guarantees.
And yet, empirically, the rotor-router final state looks much rounder that that of the
sandpile, whose boundary has clear horizontal, vertical, and slope $\pm1$ segments.

At best, this only hints at why the sandpile and rotor-router internal structures
seem to have something in common.  For now, these hints are the best I can do.

\section*{Acknowledgements}
Thanks most of all to Jim Propp, who introduced me to this lovely material, showed me
much of what appears here, and allowed and encouraged me to help spread the word.
Fond thanks to Tetsuji Miwa, whose hospitality at Kyoto University gave me the time to
think and write about it.  Thanks also to Joshua Cooper, Joel Spencer, and Matthew Cook, 
for sharing helpful comments and insights.

\end{document}